# Probability Bracket Notation for Probability Modeling

Xing M. Wang[1], Tony C. Scott[2]

## Abstract

Following the Dirac notation in Quantum Mechanics (QM), we propose the Probability Bracket Notation (PBN) by defining a probability-bra (P-bra), P-ket, P-bracket, P-identity, etc. Using the PBN, many formulae, such as normalizations and expectations in systems of one or more random variables, can now be written in abstract basis-independent expressions, which are easy to expand by inserting a proper P-identity. The time evolution of homogeneous Markov processes can also be formatted in such a way. Our system P-kets are identified with probability vectors, and our system P-bra is comparable to the Doi state function or the Peliti standard bra. In the Heisenberg picture of the PBN, a random variable becomes a stochastic process, and the Chapman–Kolmogorov equations are obtained by inserting a time-dependent P-identity. Also, some QM expressions in the Dirac notation are naturally transformed into probability expressions in PBN by a special Wick rotation. Potential applications show the usefulness of the PBN beyond the constrained domain and range of Hermitian operators on Hilbert Spaces in QM all the way to IT.

## 1. Introduction

The postulates of Quantum Mechanics (QM) were established in terms of Hilbert spaces and linear Hermitian operators. Values of physical observables such as energy and momentum were considered as eigenvalues, more precisely as spectral values of linear operators in a Hilbert space. Dirac's Vector Bracket Notation (VBN) is a very powerful tool to manipulate vectors in Hilbert spaces. It has been widely used in QM and quantum field theories. The main beauty of the VBN is that many formulae can be presented in an abstract symbolic fashion, independent of state expansions or basis selections, which, when needed, can be easily performed by inserting an appropriate v-identity operator ([1] p. 96).

$$\text{v-Bracket}: \langle \psi_A | \psi_B \rangle \implies \text{v-bra}: \langle \psi_A |, \quad \text{v-ket}: | \psi_B \rangle \tag{1}$$

$$\text{where}: \langle \psi_A | \psi_B \rangle = \langle \psi_B | \psi_A \rangle^*, \quad \langle \psi_A | = | \psi_A \rangle^\dagger \tag{2}$$

$$\text{v-basis \& v-identity}: \hat{H} | \varepsilon_i \rangle = \varepsilon_i | \varepsilon_i \rangle, \quad \langle \varepsilon_i | \varepsilon_j \rangle = \delta_{ij}, \quad \hat{I}_H = \sum_i | \varepsilon_i \rangle \langle \varepsilon_i | \tag{3}$$

$$\text{Normalization}: 1 = \langle \Psi | \Psi \rangle = \langle \Psi | \hat{I}_H | \Psi \rangle = \sum_i \langle \Psi | \varepsilon_i \rangle \langle \varepsilon_i | \Psi \rangle = \sum_i | c_i |^2 \tag{4}$$

$$\text{Expectation}: \langle H \rangle \equiv \bar{H} \equiv \langle \Psi | \hat{H} | \Psi \rangle = \sum_i \langle \Psi | \hat{H} | \varepsilon_i \rangle \langle \varepsilon_i | \Psi \rangle = \sum_i \varepsilon_i | c_i |^2 \tag{5}$$

[1] Sherman Visual Lab, Sunnyvale, CA 94085, USA; xmwang@shermanlab.com
[2] Institut für Physikalische Chemie, RWTH Aachen Univ., 52056 Aachen, Germany; tcscott@gmail.com

where † denotes the Hermitian adjoint and ∗ is the complex conjugate. However, when applying operators, the more mathematically minded worried if the transformation yielded a result contained within a Hilbert Space or subspace. They had to consider if the operator was bounded or unbounded, Referee 3: etc. One had to develop an entire spectral theory for Hermitian operators in a Hilbert space [2].

Inspired by the great success of the VBN in QM, we now propose the Probability Bracket Notation (PBN). The latter relies on a sample space, which is less constrained than a Hilbert space. Assuming X is a random variable (R.V), Ω is the set of all its outcomes, and $P(x_i)$ is the probability of $x_i \in \Omega$. Then we can make an expression like Equation (4):

$$1 = \sum_i P(x_i) = \sum_i P(x_i | \Omega) = \sum_i P(\Omega | x_i) P(x_i | \Omega) = P(\Omega | \{\sum_i |x_i) P(x_i|\} | \Omega) \leftrightarrow P(\Omega | \Omega)$$

Here, we have used the definition of conditional probability for $A, B \subseteq \Omega$ ([3], p91),

$$P(A | B) \equiv P(A \cap B) / P(B), \quad \therefore P(x_i | \Omega) = P(x_i), P(\Omega | x_i) = 1 \text{ for } \forall x_i \in \Omega \quad (6)$$

Therefore, we seem to have discovered a probability "identity operator":

$$I_X \equiv \sum_i |x_i) P(x_i| \Rightarrow P(\Omega | \Omega) = P(\Omega | I_X | \Omega) = \sum_i P(\Omega | x_i) P(x_i | \Omega) \underset{(1.7)}{=} \sum_i P(x_i) = 1$$

Then, following Dirac's notation, we define the probability bra (P-bra), P-ket, P-bracket (as conditional probabilities by nature), P-basis, the system P-ket, P-identity, normalization, expectation and more, similar but not identical to their counterparts in Equations (1) - (5). In Section 2, for systems of one R.V, we show that the PBN has an advantage similar to that of the VBN: miscellaneous probability expressions [3–5] now can be presented in an abstract way, independent of P-basis, and can be expanded by inserting a suitable P-identity.

Next, in Section 3, we investigate the time evolution of homogeneous Markov chains (HMC) [3–5]. We realize that the time evolution of a continuous time HMC can be written in a symbolic abstract expression (in Section 3.2), just like the stationary Schrödinger equation in the VBN:

$$i\hbar \frac{\partial}{\partial t} | \Psi(t)\rangle = \hat{H} | \Psi(t)\rangle = \left[ -\frac{\hbar^2}{2m} \frac{\partial^2}{\partial x^2} - V(x) \right] | \Psi(t)\rangle,$$

$$| \Psi(t)\rangle = U(t) | \Psi(0)\rangle = \exp\left[ \frac{-i}{\hbar} \hat{H} t \right] | \Psi(0)\rangle \quad (7)$$

We also find that our time dependent system P-kets can be identified with probability vectors ([3]). Our system P-bra is closely related to the state function or standard bra introduced in Doi-Peliti Techniques [6–8]. We show that by transforming from the Schrödinger picture to the Heisenberg picture in the PBN, the time dependence of a system P-ket relocates to the R.V, which becomes a stochastic process; the Chapman-Kolmogorov Theorem [4,5,9] for transition probabilities can be derived by just inserting a time dependent P-identity. Section 5 shows that a Schrödinger equation transforms to a master equation by making a special Wick rotation. Section 6 showcases the potential applications of the PBN, such as handling non-Hermitian operators and clustering text datasets. Discussion and concluding remarks are made at the end.

**2. Probability Bracket Notation and Random Variable (R.V)**
*2.1. Discrete random variable*

We define a probability space $(\Omega, X, P)$ of a discrete random variable ($R.V$, or observable) $X$ as follows: the set of all elementary events $\omega$, associated with a discrete random variable $X$, is the sample space $\Omega$, and

$$\text{For } \forall \omega_i \in \Omega, X(\omega_i) = x_i \in \Re, \quad P: \omega_i \mapsto P(\omega_i) = m(\omega_i) \geq 0, \sum_i m(\omega_i) = 1 \tag{8}$$

**Definition 1.** (*Probability event-bra and evidence-ket*): Let $A \subseteq \Omega$ and $B \subseteq \Omega$,

1. The symbol $P(A|$ represents a probability event bra, or *P*-bra;
2. The symbol $|B)$ represents a probability evidence ket, or *P*-ket.

**Definition 2.** (*Probability Event-Evidence Bracket*): The *conditional probability* (CP) of event $A$ given evidence $B$ in the sample space $\Omega$ can be treated as a *P-bracket*, and it can be split into a *P*-bra and a *P*-ket, similar to a Dirac bracket. For $A, B \subseteq \Omega$, we define:

$$P\text{-bracket: } P(A|B) \Rightarrow P\text{-bra: } P(A|, \quad P\text{-ket: }|B); \quad \text{Note: } P(A|\neq|A)^\dagger \quad (2.2a)$$

$$\text{Here: } P(A|B) \equiv \frac{P(A \cap B)}{P(B)} \quad \text{(a conditional probability by nature)} \tag{9}$$

As a CP, the *P*-bracket has the following properties for $A, B \subseteq \Omega$:

$$P(A|B) = 1 \quad if \ A \supseteq B \supset \emptyset \tag{10}$$
$$P(A|B) = 0 \quad if \ A \cap B = \emptyset \tag{11}$$

**Definition 3. System P-ket:** For any subset $E \subseteq \Omega$, the probability $P(E)$ can be written as a conditional probability, or a *P*-bracket:

$$P(E) = P(E|\Omega) \tag{12}$$

Here $|\Omega)$ is called the ***system P-ket***. The *P*-bracket defined in (9) now becomes:

$$P(A|B) = \frac{P(A \cap B)}{P(B)} = \frac{P(A \cap B|\Omega)}{P(B|\Omega)} \tag{13}$$

We have the following important property expressed in *PBN*:

$$\text{For } \forall B \subseteq \Omega \text{ and } B \neq \emptyset, \quad P(\Omega | B) = 1 \tag{14}$$

The *Bayes formula* (see [3], Sec. 2.1) now can be expressed as:

$$P(A | B) \equiv \frac{P(B | A) P(A)}{P(B)} = \frac{P(B | A) P(A | \Omega)}{P(B | \Omega)} \tag{15}$$

The set of all elementary events in Ω forms a complete, mutually disjoint set:

$$\bigcup_{\omega_i \in \Omega} \omega_i = \Omega, \quad \omega_i \cap \omega_j = \delta_{ij} \omega_i, \quad \sum_i m(\omega_i) = 1 \tag{16}$$

**Definition 4** (*Discrete P-Basis and P-Identity*): Using Eq. (8-11), we have the following properties for *basis* elements in (Ω, X, P):

$$X(\omega_j) = x_j \rightarrow X | \omega_j) = x_j | \omega_j), \quad P(\Omega | \omega_j) = 1, \quad P(\omega_i | \Omega) = m(\omega_i) \tag{17}$$

In view of the one-to-one correlation between $x_j$ and $\omega_j$, from now on, we will use $x_j$ to label basis elements, just like labeling eigenstates in Eq. (3) in VBN for QM:

$$X(x_j) = x_j \rightarrow X | x_j) = x_j | x_j), \quad P(\Omega | x_j) = 1, \quad P(x_i | \Omega) = m(x_i) \tag{18}$$

Here, *X* behaves like a right-acting operator. The complete mutually disjoint events in (16-18) form a *probability basis* (or *P-basis*) and a *P-identity*, similar to Eq. (3) in QM:

$$P(x_i | x_k) = \delta_{ik}, \quad \sum_{x \in \Omega} | x) P(x | = \sum_i | x_i) P(x_i | = I_X. \tag{19}$$

The system *P*-ket, |Ω), now can be expanded from the left as:

$$| \Omega) = I_X | \Omega) = \sum_i | x_i) P(x_i | \Omega) = \sum_i m(x_i) | x_i) \tag{20}$$

While the *system P-bra*, $P(\Omega |$, has its expansion from the right as:

$$P(\Omega | = P(\Omega | I_X = \sum_i P(\Omega | x_i) P(x_i | \underset{(14)}{=} \sum_i P(x_i | \tag{21}$$

The two expansions are quite different, and $P(\Omega | \neq [| \Omega)]^\dagger$. But their *P*-bracket is consistent with the requirement of normalization, similar to Eq. (4) in VBN:

$$1 = P(\Omega) \equiv P(\Omega | \Omega) = \sum_{i,j=1}^N P(x_i | m(x_j) | x_j) = \sum_{i,j=1}^N m(x_j) \delta_{ij} = \sum_{i=1}^N m(x_i) \tag{22}$$

**Definition 5** (*Expectation Value*): Analogous to Eq. (5) in QM, the expected value of the *R.V* or observable *X* in (Ω, X, P) now can be expressed as:

$$\langle X \rangle \equiv \bar{X} \equiv E[X] = P(\Omega \mid X \mid \Omega) = \sum_{x \in \Omega} P(\Omega \mid X \mid x) P(x \mid \Omega) = \sum_{x \in \Omega} x \, m(x) \qquad (23)$$

If $f(X)$ is a continuous function of observable $X$, then it is easy to show that:

$$\langle f(X) \rangle \equiv E[f(X)] \equiv P(\Omega \mid f(X) \mid \Omega) = \sum_{x \in \Omega} f(x) \, m(x) \qquad (24)$$

*2.2. Independent random variables:*

Let $\vec{X} = \{X_1, X_2, \ldots X_n\}$ be a vector of *independent* random variables, and the sample space (i.e., the set of possible outcomes) of $X_i$ is the set $\Omega_i$. Then, the *joint probability distribution* can be denoted as:

$$P(x_1, \ldots, x_N \mid \Omega) = P(X_1 = x_1, \ldots, X_n = x_n \mid \Omega), \quad |\Omega\rangle = |\Omega_1 \otimes \ldots \otimes \Omega_n\rangle \qquad (25)$$

The joint probability of *independent R.V* is factorable, e.g.:

$$P(x_i, x_k \mid \Omega) = P(x_i \mid \Omega_i) P(x_k \mid \Omega_k) = P(x_i) P(x_k) \qquad (26)$$

The factorable system has the following factorable expectations, for example:

$$P(\Omega \mid X_i X_j X_k \mid \Omega) = P(\Omega_i \mid X_i \mid \Omega_i) P(\Omega_j \mid X_j \mid \Omega_j) P(\Omega_k \mid X_k \mid \Omega_k) \qquad (27)$$

As an example, in *Fock space*, we have the following basis from the *occupation numbers*

$$N_i \mid \vec{n} \rangle = n_i \mid \vec{n} \rangle, \, P(\vec{n} \mid \vec{n}') = \delta_{\vec{n},\vec{n}'} = \prod_i \delta_{n_i, n'_i}, \, I_{\vec{n}} = \sum_{\vec{n}} \mid \vec{n} \rangle P(\vec{n} \mid = \prod_i \sum_{n_i} \mid n_i \rangle P(n_i \mid \qquad (28)$$

The expectation value of an occupation number now is given by:

$$\langle N_i \rangle \equiv P(\Omega \mid N_i \mid \Omega) = P(\Omega_i \mid N_i \mid \Omega_i) = P(\Omega \mid N_i I_i \mid \Omega_i) = \sum_{n_i} n_i \, P(n_i \mid \Omega_i) \qquad (29)$$

Moreover, if sets *A* and *B* are mutually independent in $\Omega$, we have the following equivalence:

$$P(A \mid B) = P(A \mid \Omega) \iff P(A \cap B \mid \Omega) = P(A \mid \Omega) P(B \mid \Omega) \qquad (30)$$

*2.3. Continuous P-basis and P-Identity*

Eq. (17-19) can be extended to probability space $(\Omega, X, P)$ of a *continuous* random variable *X*,

$$X(x) = x \to X \mid x \rangle = x \mid x \rangle, \quad P(\Omega \mid x) = 1, \quad P : x \mapsto P(x) \equiv P(x \mid \Omega) \qquad (31)$$

$$P(x \mid x') = \delta(x - x'), \quad \int_{x \in \Omega} \mid x \rangle \, dx \, P(x \mid = I_X \qquad (32)$$

We see that it is consistent with the normalization requirement:

$$P(\Omega|\Omega) = P(\Omega|I_X|\Omega) = \int P(\Omega|x)\,dx\,P(x|\Omega) = \int_{x\in\Omega} dx\,P(x|\Omega) = \int_{x\in\Omega} dx\,P(x) = 1 \quad (33)$$

The expected value $E(X)$ can be easily extended from (23):

$$\langle X \rangle \equiv \bar{X} \equiv E[X] = P(\Omega|X|\Omega) = \int_{x\in\Omega} P(\Omega|X|x)\,dx\,P(x|\Omega) = \int_{x\in\Omega} dx\,x\,P(x) \quad (34)$$

The *basis-independent expressions* in the PBN are similar to those in the Dirac VBN, and among them are the expectation and normalization formulas in the two notations:

$$PBN: \quad \langle f(X) \rangle \equiv E[f(X)] = P(\Omega|f(X)|\Omega), \quad P(\Omega|\Omega) = 1 \quad (35)$$

$$VBN: \quad \langle f(\hat{X}) \rangle \equiv E[f(\hat{X})] = \langle \psi | f(\hat{X}) | \psi \rangle, \quad \langle \psi | \psi \rangle = 1 \quad (36)$$

Of course, there exist differences between the PBN and Dirac VBN. For example, the expansions of bra, ket, and normalization have their own expressions:

$$PBN: \quad |\Omega) = \int dx|x)P(x|\Omega), \quad P(\Omega| = \int dx P(x|, \quad P(\Omega|\Omega) = \int dx P(x|\Omega) = 1 \quad (37)$$

$$VBN: \quad |\psi\rangle = \int dx|x\rangle\langle x|\psi\rangle, \quad \langle\psi| = \int dx\langle\psi|x\rangle\langle x|, \quad \langle\psi|\psi\rangle = \int dx|\langle\psi|x\rangle|^2 = 1 \quad (38)$$

*2.4. Conditional probability and expectation*

The conditional expectation of $X$ given $H \subset \Omega$ in the continuous basis (31) can be expressed as ([5], p.61):

$$E[X|H] \equiv P(\Omega|X|H) = \int P(\Omega|X|x)\,dx\,P(x|H) = \int x\,dx\,P(x|H) \quad (39)$$

$$\text{where: } P(x|H) \underset{(9)}{=} \frac{P(x \cap H|\Omega)}{P(H|\Omega)} \quad (conditional\ probability) \quad (40)$$

To become familiar with PBN, let us see two simple examples.

2.4.1. Example: Rolling a Die, (Ref. [4], Example 2.6-2.8)

A die is rolled once. We let $X$ denote the outcome of this experiment. Then, the sample space for this experiment is the 6-element set $\Omega = \{1, 2, 3, 4, 5, 6\}$. We assumed that the die was fair, and we chose the distribution function defined by $m(i) = 1/6$, for $i = 1, \ldots, 6$. Using the PBN, we have the P-identity for this sample space and the probability of each outcome:

$$\sum_{i=1}^{6} |i)P(i| = I_D \quad (41)$$

$$1 = P(\Omega|\Omega) = \sum_{i=1}^{6} P(\Omega|i)P(i|\Omega) = \sum_{i=1}^{6} P(i|\Omega) = 6p \quad (42)$$

Hence the probability for each outcome has the same value:

$$P(i) \equiv P(i|\Omega) = p = 1/6 \quad (43)$$

Its expectation value can be readily calculated:

$$P(\Omega | X | \Omega) = \sum_{i=1}^{6} P(\Omega | X | i) = \sum_{i=1}^{6} P(\Omega | i) i / 6 = \sum_{i=1}^{6} i / 6 = 21/6 = 7/2 \quad (44)$$

and the variance can be calculated as:

$$\sigma^2 = \langle X^2 \rangle - \bar{X}^2 = \sum_{i=1}^{6} i^2 / 6 - 49/4 = 91/6 - 49/4 = 45/4 \quad (45)$$

2.4.2. Example: Rolling a Die (Examples 2.1 continued, [4] 2.8 continued])
If $E$ is the event that the result of the roll is an even number, then $E = \{2, 4, 6\}$ and $P(E) = m(2) + m(4) + m(6) = 1/6 + 1/6 + 1/6 = \frac{1}{2}$. Using the PBN, the probability of event $E$ can be easily calculated as:

$$P(E) \equiv P(E | \Omega) = P(E | \hat{I} | \Omega) = \sum_{i \in E} P(E | i) P(i | \Omega) = \sum_{i \in E} P(i | \Omega) = 3p = 3/6 = 1/2 \quad (46)$$

Applying Eq. (40), we can calculate the conditional probability $P(i | E)$ as follows:

$$P(i | E) = \frac{P(E \cap i | \Omega)}{P(E | \Omega)} = \begin{cases} P(i | \Omega) / (1/2) = (1/6) / (1/2) = 1/3, & (i \text{ even}) \\ 0, & (i \text{ odd}) \end{cases} \quad (47)$$

Our discussions above can be easily extended to systems of multiple $R.V$. For example, we can introduce the following system of three independent discrete $R.V$:

$$X | x, y, z) = x | x, y, z), \quad Y | x, y, z) = y | x, y, z), \quad Z | x, y, z) = z | x, y, z) \quad (48)$$
$$P(x, y, z) = P(x, y, z | \Omega); P(\Omega | \Omega) = 1, P(\Omega | H) = 1, \text{ if } H \subseteq \Omega, \rightarrow P(\Omega | x, y, z) = 1$$

Orthonormality: $P(x, y, z' | x', y', z') = \delta_{xx'} \delta_{yy'} \delta_{zz'}$ (49)

Completeness: $\sum_{x,y,z \in \Omega} | x, y, z) P(x, y, z | = I_{X,Y,Z}$ (50)

## 3. Probability Vectors and Homogeneous Markov Chains (HMCs)

For simplicity, from now on, we will only discuss time evolution for one $R.V$. We assume our probability space $(\Omega, N, P)$ has the following stationary discrete $P$-basis from a random variable $\hat{N}$ (possible for state-labeling or occupation-number counting):

$$\hat{N} | i) = i | i), \quad P(i | j) = \delta_{ij}, \quad \sum_{i=1}^{N} | i) P(i | = I \quad (51)$$

*3.1. Discrete-time HMC*
The transition *matrix element* $p_{ij}$ is defined as ([3], p.407):

$$p_{ij} \equiv P(N_{t+1} = j | N_t = i) \equiv P(j, t+1 | i, t), \quad \sum_{j=1}^{N} p_{ij} = 1 \quad (52)$$

In matrix form, if we define a *probability row vector* (PRV) at t = 0 as $\langle u^{(0)}|$, then matrix $P$ acting on it from the right for $k$ times gives the PRV at time = $k$ ([3], theorem 11.2):

$$\langle u^{(k)}| = \langle u^{(0)}|P^k, \text{ or}: \quad \langle u^{(k)}|i\rangle = u^{(k)}{}_i = u^{(0)}{}_j P^k{}_{ji} \tag{53}$$

**Definition 6.** *Time-dependent System P-ket*: we use the following system *P*-ket, to represent a *probability column vector* in probability space $(\Omega, X, P)$:

$$|\Omega_t\rangle = \sum_i^N |i\rangle P(i|\Omega_t) = \sum_i^N m(i,t)|i\rangle, \quad P(\Omega|\Omega_t) = \sum_i^N m(i,t) = 1 \tag{54}$$

The time evolution equation (53) can now be written in a basis-independent way:

$$|\Omega_t\rangle = (P^T)^t |\Omega_0\rangle \equiv \hat{U}(t)|\Omega_0\rangle, \quad t \in \mathbb{N}; \quad |\Omega_0\rangle = [\langle u^{(0)}|]^T \tag{55}$$

**Definition 7** (*Time-dependent Expectation*): The expectation value of a continuous function $f$ of occupation number $\hat{N}$ in $(\Omega, N, P)$ can be expressed as:

$$\langle f(\hat{N})\rangle = P(\Omega|f(\hat{N})|\Omega_t) = \sum_i P(\Omega|f(i)|i)P(i|\Omega_t) = \sum_i f(i)m(i,t) \tag{56}$$

We can map the *P*-bra and *P*-ket into Hilbert space by using Dirac's notation:

$$P(\Omega| = \sum_i P(i| \leftrightarrow \langle\Psi| = \sum_i \langle i|, \quad |\Omega_t\rangle \leftrightarrow |\Psi_t\rangle = \sum_i |i\rangle\langle i|\Psi_t\rangle = \sum_i c(i,t)|i\rangle \tag{57}$$

Then, the expectation Eq. (56) can be rewritten in Dirac's notation as:

$$\langle f(\hat{N})\rangle = \langle\Psi_t|f(\hat{N})|\Psi_t\rangle = \sum_i \langle\Psi_t|f(i)|i\rangle\langle i|\Psi_t\rangle = \sum_i f(i)|c(i,t)|^2 \tag{58}$$

$$\therefore P(i|\Omega_t) \equiv m(i,t) = |c(i,t)|^2 \equiv |\langle i|\Psi_t\rangle|^2, \quad t \in \mathbb{N} \tag{59}$$

One can see in Equation (59) that the local phases of quantum states are not contained in the probability distribution, and they have no contribution to the distribution:

$$VBN : |\Psi_t\rangle = \sum_i c_i(t)|\psi_i\rangle \underset{\text{one to many}}{\overset{\text{many to one}}{\rightleftarrows}} PBN : |\Omega_t\rangle = \sum_i P(i,t)|i\rangle; P(i,t) = |c_i(t)|^2$$

This means a certain loss of quantum information, and it is well known that these local phase factors are crucial in explaining experiments like the quantum interference [10] and the Aharonov–Bohm effect [11].

*3.2. Continuous-time HMC*

The time-evolution equation of a continuous-time HMC with a discrete-basis (see Eq. (85) or [5], p.221) can be written as:

$$\frac{\partial}{\partial t} p_j(t) = \sum_k p_k(t) q_{kj} \quad \Rightarrow \quad \frac{\partial}{\partial t} P(j|\Omega_t) = \sum_k \hat{Q}^T_{jk} P(k|\Omega_t) \quad (60)$$

$$= \sum_k P(j|\hat{Q}^T|k) P(k|\Omega_t) = P(j|\hat{Q}^T I|\Omega_t) = P(j|\hat{Q}^T|\Omega_t) \equiv P(j|\hat{L}|\Omega_t) \quad (61)$$

Eq. (60-61) leads to a basis-independent **master equation** (see [8] Eq. (2.12)):

$$\frac{\partial}{\partial t}|\Omega_t) = \hat{L}|\Omega_t), \quad |\Omega_t) = \hat{U}(t)|\Omega_0) = e^{\hat{L}t}|\Omega_0), \quad t \geq 0 \quad (62)$$

It looks just like Schrodinger's equation (7) of a conserved quantum system in Dirac's notation. With the discrete *P*-basis of Fock space in Eq. (28), Eq. (56-59) can now be written as:

$$|\Omega_t) = \sum_{\vec{n}} m(\vec{n}, t)|\vec{n}), \quad P(\Omega| = \sum_{\vec{n}} (\vec{n}|, \quad \therefore \langle f(\vec{n}) \rangle = P(\Omega|f(\vec{n})|\Omega_t) \quad (63)$$

Doi's definition of the state function $\langle s |$ and state vector $| F(t) \rangle$ ([8, 9]) corresponds to our system *P*-bra and *P*-ket, respectively:

$$P(\Omega| = \sum_{\vec{n}} P(\vec{n}| \quad \leftrightarrow \quad \langle s | \equiv \sum_{\vec{n}} \langle \vec{n} |$$

$$|F(t)\rangle \equiv \sum_{\vec{n}} P(\vec{n},t)|\vec{n}\rangle \quad \leftrightarrow \quad |\Omega_t) = \sum_{\vec{n}} |\vec{n}) P(\vec{n}|\Omega_t)$$

$$\therefore \langle \hat{B}(\vec{n}) \rangle = \langle s | \hat{B}(\vec{n}) | F(t) \rangle = P(\Omega|\hat{B}(\vec{n})|\Omega_t) \quad (64)$$

Note that the vector-basis here corresponds to the *P*-basis in Eq. (28):

$$\hat{n}_i|\vec{n}\rangle = n_i|\vec{n}\rangle, \quad \sum_{\vec{n}} |\vec{n}\rangle\langle\vec{n}| = I_{\vec{n}}, \quad \langle\vec{n}|\vec{n}'\rangle = \delta_{\vec{n},\vec{n}'} = \prod_{i=1} \delta_{n_i,n_i'} \quad (65)$$

In Peliti's formalism ([8] p.1472), the vector-basis (from population operator *n*) is normalized in a special way, so the expansion of the system P-bra is also changed:

$$\sum_n |n\rangle \frac{1}{n!} \langle n| = I_n, \quad \langle m|n\rangle = n! \delta_{m,n} \quad (66)$$

$$P(\Omega| = P(\Omega|I_n = \sum_n^{\infty} P(\Omega|n) \frac{1}{n!} P(n| \underset{(2.3)}{=} \sum_n^{\infty} P(n| \frac{1}{n!} \quad (67)$$

Eq. (67) can be identified with the *standard bra*, introduced in [8] Eq. (2.27-28):

$$P(\Omega| = \sum_n \frac{1}{n!} P(n| \quad \leftrightarrow \quad \langle| \equiv \sum_n \frac{1}{n!} \langle n|$$

$$|\Omega) = \sum_n |n) \frac{1}{n!} P(n|\Omega) \quad \leftrightarrow \quad |\phi\rangle = \sum_m |m\rangle \frac{1}{m!} \langle m|\phi\rangle$$

$$\therefore E[\hat{A}] \equiv \langle \hat{A} \rangle = \langle|\hat{A}|\phi\rangle = P(\Omega|\hat{A}|\Omega) \quad (68)$$

*3.3. The Heisenberg Picture*

We call Eq. (55) and (62) the evolution equations in the *Schrodinger picture*. Now we introduce the *Heisenberg picture* of an *R.V* (or observable), similar to what is used in QM ([1] p.541):

$$|\Omega_t) = \hat{U}(t)|\Omega_0) \quad \Rightarrow \quad \hat{X}(t) = \hat{U}^{-1}(t)\hat{X}\hat{U}(t) \tag{69}$$

Based on $\hat{U}(t)$, we can introduce the following time-dependent P-basis:

$$|x,t) = \hat{U}^{-1}(t)|x), \quad P(x|\hat{U}(t) = P(x,t|, \quad P(x,t|x',t) = P(x|x'), \quad P(\Omega|x,t) = 1 \tag{70}$$

$$P(x',t|\hat{X}(t)|x,t) = P(x'|\hat{U}(t)\hat{U}^{-1}(t)\hat{X}\hat{U}(t)\hat{U}^{-1}(t)|x) = P(x'|\hat{X}|x) = xP(x'|x) \tag{71}$$

The probability density can now be interpreted in the two pictures:

$$P(x,t) \equiv P(x|\Omega_t) = P(x|\hat{U}(t)|\Omega_0) = P(x,t|\Omega_0) = P(x,t|\Omega) \tag{72}$$

In the last step, we have used the fact that in the Heisenberg picture, $|\Omega_0) = |\Omega)$.

**Definition 8.** *The Time-dependent P-Identity*: Eq. (69-71) also provides us with a time-dependent *P*-identity in the Heisenberg Picture:

Discrete: $I_X(t) = \hat{U}^{-1}(t)I_X\hat{U}(t) = \hat{U}^{-1}(t)\sum_i |x_i)P(x_i|\hat{U}(t) \underset{(3.18)}{=} \sum_i |x_i,t)P(x_i,t|$

Continuous; $I_X(t) = \hat{U}^{-1}(t)I_X\hat{U}(t) = \hat{U}^{-1}(t)\int |x)\,dx\,P(x|\hat{U}(t) = \int |x,t)\,dx\,P(x,t| \tag{73}$

Now, the expectation value of the stochastic process $X(t)$ can be expressed as:

$$P(\Omega|\hat{X}(t)|\Omega) = P(\Omega|\hat{X}(t)I_X(t)|\Omega) = \int dx\,P(\Omega|\hat{X}(t)|x,t)P(x,t|\Omega)$$
$$= \int dx\,x\,P(x,t|\Omega) = \int dx\,x\,P(x|\Omega_t) = P(\Omega|X|\Omega_t) \tag{74}$$

This suggests that a Markov *stochastic process* $X(t)$ can be thought of as an *operator in the Heisenberg picture*, and its expectation value can be found from its Schrodinger picture. Additionally, if a Markov process $X(t) \equiv X_t$ is homogeneous, with independent and stationary increments ([4], p15), we can always set $X_0 = 0$, and obtain the following useful property:

$$P(X_{t+s} - X_s = x) = P(X_t - X_0 = x|\Omega)$$
$$= P(X_t = x|\Omega) = P(x,t|\Omega) = P(x|\Omega_t) \equiv P(x,t) \tag{75}$$

Moreover, there is a relation between their transition probability and probability density:

$$P(x_2, t_1 | x_2, t_1) = P(x_2 - x_1, t_2 - t_1 | \Omega) \equiv P(x_2 - x_1, t_2 - t_1), \quad (t_1 < t_2) \qquad (76)$$

*3.4 Chapman-Kolmogorov Equations* ([5], [9])

They can now be easily obtained by inserting our time-dependent *P*-identity in Eq. (73). For an HMC of *discrete R.V* ([5] p.174):

$$p^{m+n}{}_{ij} \equiv P(j, m+n | i, 0) = P(j, m+n | \hat{I}(m) | i, 0) = \sum_k P(j, m+n | k, m) P(k, m | i, 0)$$
$$\underset{\text{HMC}}{=} \sum_k P(j, n | k, 0) P(k, m | i, 0) \underset{(4.2)}{=} \sum_k p^m{}_{ik} \, p^n{}_{kj} \quad \text{(Discrete time)} \qquad (77)$$

$$p_{ij}(t+s) \equiv P(j, t+s | i, 0) = P(j, t+s | \hat{I}(s) | i, 0) = \sum_k P(j, t+s | k, s) P(k, s | 0, i)$$
$$\underset{\text{HMC}}{=} \sum_k P(j, t | k, 0) P(k, s | 0, i) = \sum_k p_{ik}(s) \, p_{kj}(t) \quad \text{(Continuous time)} \qquad (78)$$

For a Markov process of both *continuous R.V* and time [9] (Eq. (3.9), p.31):

$$P(x, t | y, s) = P(x, t | \hat{I}(\tau) | y, s) = \int P(x, t | z, \tau) dz \, P(z, \tau | y, s) \quad \text{where } t > \tau > s \qquad (79)$$

*Absolute probability distribution* (APD)

Likewise, we can find the time evolution of APD for an HMC simply by inserting the P-identity $I(t=0) = I(0)$: For an HMC with discrete states ([5] pp.174, 214):

$$\textit{Discrete-time}: \qquad P(i, m) = P(i, m | \Omega) = P_i^{(m)}, \quad P(i, 0) \equiv p_i^{(0)}$$
$$\therefore P_i^{(m)} = P(i, m | I(0) | \Omega) = \sum_k P(i, m | k, 0) P(k, 0 | \Omega) = \sum_k p_k^{(0)} p^m{}_{ki} \qquad (80)$$

Note that Eq. (80) is identical to Eq. (53):

$$\textit{Continuous time}: \; P(x_i, t) = P(x_i, t | \Omega) = P(x_i | \Omega_t), \quad P(x_i, 0) \equiv p_i(0)$$
$$\therefore P(x_i, t) = P(x_i, t | I(0) | \Omega) = \sum_k P(x_i, t | x_k, 0) P(x_k, 0 | \Omega) = \sum_k p_k(0) p_{ki}(t) \qquad (81)$$

For a Markov process of both *continuous R.V* and time ([9], p36):

$$P(x, t) = P(x, t | \Omega) = P(x | \Omega_t), \quad P(x, 0) \equiv p^{(0)}(x)$$
$$\therefore P(x, t) = P(x, t | I(0) | \Omega) = \int dy P(x, t | y, 0) P(y, 0 | \Omega) = \int dy P(x, t | y, 0) p^{(0)}(y) \qquad (82)$$

*3.5. Kolmogorov Forward and Backward Equations*

If the Markov chains are *stochastically continuous*, then for infinitesimal $h$, the transition probability has the Tailor expansions ([4] Sec. 6.8; [5] p.217):

$$p_{ij}(h) = p_{ij}(0) + p_{ij}'(0)h + o(h^2) = \delta_{ij} + q_{ij}h + o(h^2) \qquad (83)$$

Then, using the Chapman-Kolmogorov Eq. (78), we have:

$$p_{ij}(t+h) = \sum_k p_{ik}(t) \, p_{kj}(h) = \sum_k p_{ik}(t) \, (\delta_{kj} + q_{kj}h + o(h^2))$$
$$= p_{ij}(t) + \sum_k p_{ik}(t) \, (q_{kj}h + o(h^2)) \tag{84}$$

Therefore, we get the following Kolmogorov *Forward* equations:

$$p'_{ij}(t) = \lim_{h \to 0}[(p_{ij}(t+h) - p_{ij}(t))/h] = \sum_k p_{ik}(t) \, q_{kj} \tag{85}$$

Similarly, we can derive the Kolmogorov *Backward* equations:

$$p_{ij}(h+t) = \sum_k p_{ik}(h) \, p_{kj}(t) \Rightarrow p'_{ij}(t) = \sum_k q_{ik} \, p_{kj}(t) \tag{86}$$

Using Equations (81) and (83), one can extend Equation (62) to a continuous-time HMC.

*3.6. Transition Probability and Path Integrals*

From Eq. (62) and (70), the transition probability of a Markov process can be expressed as follows (assuming $t_a < t_0 < t_b$):

$$P(x_b, t_b \mid x_a, t_b) = P(x_b \mid U(t_b, t_0)U^{-1}(t_a, t_0) \mid x_a)$$
$$= P(x_b \mid U(t_b, t_a) \mid x_a) = P(x_b \mid e^{\int_{t_a}^{t_b} dt \, L} \mid x_a) \tag{87}$$

We can divide the time interval $[t_a, t_b]$ into small pieces and insert *P*-identity *N* times:

$$\Delta t = t_{n+1} - t_n = (t_b - t_a)/(N+1) > 0, \quad t_N = t_b, t_0 = t_a, \tag{88}$$

$$P(x_b, t_b \mid x_a, t_a) = P(x_b, t_b \mid I(t_N)..I(t_1) \mid x_a, t_a) = \prod_{n=1}^{N} \int_{-\infty}^{\infty} dx_n P(x_n, t_n \mid x_{n-1}, t_{n-1})$$
$$= \prod_{n=1}^{N} \int_{-\infty}^{\infty} dx_n P(x_n \mid \hat{U}(t_n, t_{n-1}) \mid x_{n-1}) = \prod_{n=1}^{N} \int_{-\infty}^{\infty} dx_n P(x_n \mid e^{L\Delta t} \mid x_{n-1}) \tag{89}$$

It perfectly matches the starting equation of Feynman's path integral for the transition amplitude in QM ([12] Eq. (2.4), p.90):

$$\langle x_b, t_b \mid x_a, t_a \rangle = \prod_{n=1}^{N} \int_{-\infty}^{\infty} dx_n \langle x_n \mid \hat{U}(t_n, t_{n-1}) \mid x_{n-1} \rangle = \prod_{n=1}^{N} \int_{-\infty}^{\infty} dx_n \langle x_n \mid e^{-\frac{i}{\hbar}\hat{H}\Delta t} \mid x_{n-1} \rangle \tag{90}$$

Because the special Wick rotation can directly transform Eq. (89) to Eq. (90) and vice visa, it might imply possibly a new connection between the Hilbert space and the probability space for certain quantum systems:

$$\langle x_b, t_b | x_a, t_a \rangle \underset{it \leftarrow t}{\overset{it \to t}{\rightleftarrows}} P(x_b, t_b | x_a, t_a)$$

However, at the moment, this is just a potential application, and we need to investigate this further. For a *free particle*, $V(x) = 0$, the Schrodinger equation and resulting transition amplitude from Feynman's path integral are (see [12] Eq. (2.125)):

$$i\hbar \frac{\partial}{\partial t} \langle x | \Psi(t) \rangle = \langle x | \hat{H} | \Psi(t) \rangle = -\frac{\hbar^2}{2m} \frac{\partial^2}{\partial x^2} \psi(x,t), \quad \langle x | \Psi(t) \rangle = \psi(x,t) \tag{91}$$

$$\langle x_b, t_b | x_a, t_a \rangle = \sqrt{\frac{m}{2\pi i\hbar(t_b - t_a)}} \exp\left[\frac{im(x_b - x_a)^2}{2\hbar(t_b - t_a)}\right] \tag{92}$$

### 4. Examples of Homogeneous Markov Processes

Now, let us observe some important examples of the HMC using the PBN.

*Poisson Process* ([4] p.250; [5] p.161)

This is a counting process, $N(t)$, having the following properties:

(1). $\{N(t), t \geq 0\}$ is a non-negative process with independent increments and $N(0) = 0$;
(2). It is homogeneous, and its probability distribution is given by:

$$m(k,t) \equiv P([N(t+s) - N(s) = k] | \Omega) \underset{\substack{independent \\ increament}}{=} P([N(t) - N(0) = k] | \Omega)$$

$$\underset{N(0)=0}{=} P([N(t) = k] | \Omega) \equiv P(k | \Omega(t)) \underset{\substack{Poisson \\ Distribution}}{=} \frac{(\lambda t)^k}{k!} e^{-\lambda t} \tag{93}$$

It can be shown ([5] p.161) that:

$$\mu(t) \equiv \bar{N}(t) = \sum_i k\, m(k,t) = \lambda t; \quad \sigma^2(t) \equiv P(\Omega | [N(t) - \bar{N}(t)]^2 | \Omega) = \lambda t \tag{94}$$

It can be shown (see [5] p.215; [13] Theorem 1.5, p.6) that the Poisson Process has Markov property, and its transition probability is:

$$P([N(t+s) = j] | N(t) = i) = P([N(t+s) - N(t) = j - i] | \Omega) = \frac{(\lambda t)^{j-i}}{(j-i)!} e^{-\lambda t}, \text{ if } j \geq i$$

$$p_{ij}(t) = 0, \text{ if } j < i \tag{95}$$

*Wiener-Levy Process* (see [5] p.159; [9] §3.6, p.32)

This is a homogeneous process $\{W(t), t \geq 0\}$ with *independent and stationary increments* and $W(0) = 0$. Its probability density is a normal (Gaussian) distribution $N(0, t\sigma^2)$:

$$P(x,t) \equiv P([W(t+s)-W(s)=x]|\Omega) \underset{hom ogeneous}{=} P([W(t)-W(0)=x]|\Omega)$$

$$\underset{X(0)=0}{=} P([W(t)=x]|\Omega) \equiv P(x,t|\Omega) \equiv P(x|\Omega(t)) \underset{\substack{Normal\\Distribution}}{=} \frac{1}{\sqrt{2\pi t}\,\sigma}\exp\left[-\frac{x^2}{2t\sigma^2}\right] \quad (96)$$

Its stationary increment is defined as:

$$P(x_2,t_2 | x_1,t_1) = \frac{1}{\sqrt{2\pi(t_2-t_1)}\,\sigma}\exp\left[-\frac{(y_2-y_1)^2}{2(t_2-t_1)\sigma^2}\right], \quad (t_1 < t_2) \quad (97)$$

We see that the Wiener-Levy process satisfies Eq. (76).

$$P(x_2,t_2 | x_1,t_1) = P(x_2-x_1, t_2-t_1), \quad (t_1 < t_2) \quad (98)$$

It can be verified ([5] p.62) that it is an $N(0,\sigma^2)$:

$$\tilde{\mu}(t) \equiv P(\Omega|W(t)|\Omega) = 0, \quad \tilde{\sigma}^2(t) \equiv P(\Omega|W(t)^2|\Omega) = \sigma^2 \quad (99)$$

*Brownian motion* ([4] Sec. 10.1, p.524; [9] p.6, p.42)

The stochastic process $X(t)$ has $X(0) = 0$, has *stationary and independent increments* for $t \geq 0$ and its density function for $t > 0$ is a normal distribution $N(t\mu, t\sigma^2)$, or:

$$P(x,t) \equiv P(x|\Omega(t)) = \frac{1}{\sqrt{2\pi t}\,\sigma}\exp\left[-\frac{(x-\mu t)^2}{2t\sigma^2}\right] \quad (100)$$

Brownian motions have *stationary and independent* increments ([4] Sec. 10, pp. 524, 529), so they are HMC and satisfy Eq. (76). They are the solution of the following master equation (Einstein's diffusion equation, see [9] p.6) with a drift speed $\mu$:

$$\frac{\partial}{\partial t}P(x,t) = D\frac{\partial^2}{\partial x^2}P(x,t), \quad D = \frac{\sigma}{\sqrt{2}}. \quad (101)$$

Here, the constant $D$ is called the diffusion coefficient. Eq. (101) can be interpreted as a special HMC case of Eq. (62), given in the *x*-basis:

$$\frac{\partial}{\partial t}P(x|\Omega_t) = \int dx' P(x|\hat{L}|x')P(x'|\Omega_t), \quad P(x|\hat{L}|x') = D\frac{\partial^2}{\partial x^2}\delta(x-x') \quad (102)$$

Eq. (102) closely resembling Schrodinger's Eq. (7) for a free particle in the *x*-basis:

$$i\hbar\frac{\partial}{\partial t}\langle x|\Psi(t)\rangle = \int dx'\langle x|\hat{H}|x'\rangle\langle x'|\Psi(t)\rangle, \quad P(x|\hat{H}|x') = \frac{-\hbar^2}{2m}\frac{\partial^2}{\partial x^2}\delta(x-x') \quad (103)$$

## 5. Special Wick Rotation, Time Evolution and Induced Diffusions

The *induced microscopic diffusion* is defined by the following equation:

$$\frac{\partial}{\partial t}P(x,t) = \hat{G}(x,\hat{\kappa})P(x,t), \quad \hat{G} = \frac{1}{2\mu_\hbar}\frac{\partial^2}{\partial x^2} - \mu_\hbar u(x), \quad \mu_\hbar \equiv \frac{m}{\hbar} \qquad (104)$$

The *Special Wick Rotation* (SWR), caused by imaginary time rotation *Wick rotation* ($it \to t$) [14], is defined by:

$$\text{SWR:}\quad it \to t, \quad |\psi(t)\rangle \to |\Omega_t), \quad \langle x_b, t_b | x_a, t_a \rangle \to P(x_b, t_b | x_a, t_a) \qquad (105)$$

Under SWR, Schrodinger Eq. (7) is shifted to induced micro diffusion (104):

$$\frac{\partial}{i\partial t}|\psi(t)\rangle = \frac{-1}{\hbar}\hat{H}|\psi(t)\rangle \to \frac{\partial}{\partial t}|\Omega_t) = \hat{G}|\Omega_t) \qquad (106)$$

$$-\frac{1}{\hbar}\hat{H} = \frac{\hbar \partial_x^2}{2m} - \frac{m}{\hbar}\frac{V(x)}{m} \equiv \frac{\partial_x^2}{2\mu_\hbar} - \mu_\hbar u(x) \to \hat{G} = \frac{1}{2\mu_\hbar}\frac{\partial^2}{\partial x^2} - \mu_\hbar u(x) \qquad (107)$$

The simplest example is the induced micro *Einstein-Brown motion* when $u(x) = 0$ in Eq. (104). Applying (105) to the transition amplitude in Eq. (92), we get the transition probability:

$$P(x_b, t_b | x_a, t_a) = \sqrt{\frac{1}{4\pi D_\hbar (t_b - t_a)}} \exp\left[-\frac{(x_b - x_a)^2}{4D_\hbar (t_b - t_a)}\right] \qquad (108)$$

Here, we have introduced the *induced micro-diffusion coefficient*,

$$D_\hbar \equiv 1/(2\mu_\hbar) = \hbar/(2m) \qquad (109)$$

Similarly, applying (105) to the Schrodinger for a free particle, Eq. (103), we obtain Einstein's diffusion equation (102) with $D_\hbar \to D$.

*A Special non-Hermitian Case*

We can apply Dirac notation and the PBN together to solve the following special quantum system with a non-Hermitian Hamiltonian:

$$\hat{H} = \hat{H}_1 - i\hat{H}_2; \quad \hat{H}_1 = \hat{H}_1^\dagger; \quad \hat{H}_2 = \hat{H}_2^\dagger; \quad [\hat{H}_1, \hat{H}_2] = 0 \qquad (110)$$

Then, the time-evolution equation can be written as:

$$i\hbar\partial_t \left[|\Psi_t\rangle |\Omega_t)\right] = \left[\hat{H}_1 |\Psi_t\rangle\right]|\Omega_t) + |\Psi_t\rangle \left[-i\hat{H}_2 |\Omega_t)\right] \qquad (111)$$

It leads to two equations:

$$\partial_t |\Psi_t\rangle = -\frac{i}{\hbar}\hat{H}_1|\Psi_t\rangle, \quad \partial_t |\Omega_t\rangle = \frac{-1}{\hbar}\hat{H}_2|\Omega_t\rangle \underset{(106)}{=} \hat{G}|\Omega_t\rangle \tag{112}$$

The first is an ordinary Schrodinger equation, while the second is a mater equation for an induced micro diffusion. The product of their solutions is the solution of Eq. (111). Suppose we have $\hat{H}_1(x), \hat{H}_2(y)$; then the product of Eq. (89) and (90) gives the path-integral expression of the system:

$$\left[\langle x_b, t_b | P(y_b, t_b |\right] e^{-\int_{t_a}^{t_b} dt \frac{i}{\hbar}\left[\hat{H}_1(x) - i\hat{H}_2(y)\right]} \left[| x_a, t_a\rangle | y_a, t_a)\right]$$
$$= \prod_{m=1}^{N}\int_{-\infty}^{\infty} dx_m \langle x_m | e^{-\frac{i}{\hbar}\hat{H}_1(x)\Delta t} | x_{m-1}\rangle \prod_{n=1}^{N}\int_{-\infty}^{\infty} dy_n P(y_n | e^{\hat{G}(y)\Delta t} | y_{n-1}) \tag{113}$$

Moreover, since both $\hat{H}_1$ and $\hat{H}_2$ are Hermitian, we assume that they have complete sets of orthogonal eigenstates like:

$$\hat{H}_1|\psi_k\rangle = \varepsilon_k |\psi_k\rangle, \quad \hat{H}_2|\varphi_\mu) = \lambda_\mu |\varphi_\mu) \tag{114}$$

Then, we obtain the following solutions of Eq. (112):

$$|\psi_i(t)\rangle = e^{-i\varepsilon_k t/\hbar}|\psi_k\rangle, \quad |\varphi_\mu(t)) = e^{-\lambda_\mu t/\hbar}|\varphi_\mu) \tag{115}$$

It leads to the expectation value of $H$ given a mixed product state $|\psi_i(t)\rangle|\varphi_\mu(t))$ as:

$$E_{k,\mu}[\hat{H}] = \langle \psi_k(t)|\hat{H}_1|\psi_k(t)\rangle - iE[H_2|\varphi_{\mu,t}] = \varepsilon_k - i\lambda_\mu \tag{116}$$

Here, we have used the conditional prediction of $\hat{H}_2$ given $\varphi_\mu(t)$, which can be easily obtained from Eq. (40):

$$E[H_2|\varphi_{\mu,t}] = \sum_{\alpha=\mu}\lambda_\alpha P(\varphi_{\alpha,t})/P(\varphi_{\mu,t}) = \lambda_\mu$$

The imaginary part in Eq. (116) looks strange, but it would appear in expectation values of similar complex Hamiltonians, as in [15], p. 290 and [16], p.3, because $\hat{H}$ is not Hermitian. However, we know that $\hat{H}^\dagger\hat{H}$ is Hermitian. Therefore, it is reasonable to use the following expression:

$$\left\{E_{k,\mu}[\hat{H}^\dagger\hat{H}]\right\}^{1/2} = \left\{E_{k,\mu}[(\hat{H}_1)^2 + (\hat{H}_2)^2]\right\}^{1/2} = \left\{\varepsilon_k^2 + \lambda_\mu^2\right\}^{1/2} \tag{117}$$

It is interesting to note that if the initial condition of $|\Omega_t\rangle$ is a linear combination of a set of P-kets $|\varphi_\mu\rangle$, then because each P-ket has its probability proportional to $\exp[-\lambda\mu t/h]$, the P-ket with the lowest eigenvalue $\lambda_{min}$ will eventually dominate, and the expectation value of energy $H_2$ will approach $\lambda_{min}$. Therefore, our master Equation (104) describes the processes that will ultimately jump to its lowest energy level, included in its initial condition! Of course, further research on this is required.

## 6. Potential Applications

Hermitian operators are sufficient for "pure" eigenvalue states for closed systems where the energies are conserved and real-valued. However, for mixed states and in a number of physical circumstances, non-Hermitian operators have had to be considered. It is well known that any non-Hermitian linear operator $\hat{C}$ can be expressed in the form:

$$\hat{C} = \hat{A} + i\hat{B} \tag{118}$$

where $\hat{A}$ and $\hat{B}$ are Hermitian by letting $\hat{A} = (\hat{C} + \hat{C}^\dagger)/2$ and $\hat{B} = (\hat{C} - \hat{C}^\dagger)/2i$. Note that $\hat{B} \to i\hat{B}$ can be treated as a Wick rotation in the path-integral application of the PBN. Section 5 dovetails into the matter of general linear operators where there have been a number of applications. We mention:

1. The method of "complex scaling" applied to quantum mechanical Hamiltonians was a "hot" area in the area of atomic and molecular physics during the 1970s and 1980s and involved non-Hermitian linear operators. We cite an application by the mathematician Barry Simon on complex scaling to non-relativistic Hamiltonians for molecules [17].
2. Another application requiring such an operator was performed by Botten et al. [18]. They showed how to solve a practical problem involving wave scattering using a bi-orthogonal basis, with a VPN bra basis and a ket basis consisting of different functions. In a unitary problem, these VPN bra and ket basis functions would be the same. Here, the Helmholtz equation Laplacian $\nabla^2 f + k^2 f = 0$ has a wave number $k$, which is complex. The imaginary part of $k$ indicates loss or gain depending on its sign.
3. K. G. Zloshchastiev has performed many applications on the general density operator approach with non-Hermitian Hamiltonians of the form Equation (118) applied to, e.g., open dissipative systems, which automatically deal with mixed states (see [19–21]) and von Neumann quantum entropy [22]. For non-Hermitian, open/dissipative systems, there is also the work of Refs. [23, 24].
4. Consider a rectangular real data matrix $\hat{Q}$. Its similarity matrix (or adjacency matrix) $\hat{S}$ and corresponding row stochastic (Markov) matrix $\hat{R}$ is defined by:

$$\hat{S} = \hat{Q} * \hat{Q}^T, \quad R_{i,j} = S_{i,j} = \sum_k S_{i,k} \tag{119}$$

The symmetric matrix $\hat{s}$ has real eigenvalues, and as does $\hat{R}$, even though the latter is non-symmetric. $\hat{R}$ is a transition matrix in the language of QM but also the key operator in the Meila–Shi algorithm of spectral clustering, as well as quantum clustering with IT applications in science, engineering, (unstructured) text using a "bag-of-words" model [25–28] and even medicine [29]. Probably, the PBN may provide us with some new approaches to quantum clustering.

If the real data matrix $\hat{Q}$ in Equation (119) is huge but very sparse (as is often the case in a "bag-of-words" model [27]), we may use the PBN to find other alternative algorithms for text document clustering. Suppose that the dataset has a vocabulary of $N$ labeled keywords, which serve as the P-basis in Equation (28), and $q_{\mu,k}$ represents the frequency of the $k^{th}$ keyword in document $Q_\mu$. Since $q_{\mu,k}$ can be greater than 1, they behave like the occupation numbers of a Boson system in Quantum Field Theory. The conditional probability of finding doc $Q_\mu$ given doc $Q_\nu$ is:

$$P(Q_\mu | Q_\nu) = \sum_k^N P(Q_\mu | k) P(k | Q_\nu) \underset{(14)}{=} \sum_{k \in Q_\mu} P(k | Q_\nu); \quad P(k | Q_\mu) = \frac{q_{\mu,k}}{\sum_{k \in Q_\mu} q_{\mu,k}} \quad (120)$$

Now, for example, we can define the relevance of two docs and obtain its expression as:

$$R_{\mu,\nu} \equiv \frac{1}{2}[P(Q_\mu | Q_\nu) + P(Q_\nu | Q_\mu)] = \frac{1}{2}\left[\sum_{k \in Q_\mu} P(k | Q_\nu) + \sum_{k \in Q_o} P(k | Q_\mu)\right] \quad (121)$$

This algorithm may be very effective for huge datasets in numerous applications.

## 7. Summary and Discussion

Inspired by Dirac's notation used in quantum mechanics (QM), we proposed the Probability Bracket Notation (PBN). We demonstrated that the PBN could be a very useful tool for symbolic representation and manipulation in probability modeling, like the various normalization and expectation formulas for systems of single or multiple random variables, as well as the master equations of homogeneous Markov processes. We also show that a stationary Schrödinger Equation (103) naturally becomes a master equation under the special Wick rotation (see Equations 106-107).

We have shown the similarities between many QM expressions in VBN and related probabilistic expressions in PBN, which might provide a beneficial bridge connecting the quantum world with the classical one. We also showed how the PBN could be applied to problems that require non-Hermitian operators. We make no pretense that the PBN creates new Physics. Rather, it provides a notational interdisciplinary "umbrella" for various statistical and physical processes and thereby opens up the possibility of their synthesis and integration.

## 8. Conclusions

Of course, more investigations need to be done to verify the consistency (or correctness), usefulness, and limitations of our proposal of the PBN. Still, it is intriguing in its possibility of expressing formulations from classical statistics and QM and handling non-Hermitian operators.

**Abbreviations**

| | |
|---|---|
| APD | Absolute Probability Distribution |
| CP | Conditional Probability |
| HMC | Homogeneous Markov Chains |
| IT | Information Technology |
| PBN | Probability Bracket Notation |
| P-basis | Probability basis |
| P-bra | Probability (event) bra |
| P-identity | Probability identity |
| P-ket | Probability (event) ket |
| PRV | Probability Row Vector |
| QM | Quantum Mechanics |
| R.V | Random Variable |
| SWR | Special Wick Rotation |
| VBN (Dirac) | Vector Bracket Notation |

**References**

1. Kroemer, H. Quantum mechanics for engineering: materials science and applied physics; Pearson: Upper Saddle River, NJ, 1994.
2. Rudin, W. Functional Analysis; McGraw-Hill series in higher mathematics, McGraw-Hill: New York, NY, 1973.
3. Grinstead, C.M.; Snell, J.L., Eds. Introduction to Probability; American Mathematical Society: Providence, RI, 1997.
4. Ross, S.M. Introduction to probability models, ISE, 9 ed.; Introduction to Probability Models, Academic Press: San Diego, CA, 2006.
5. Ye, E.; Zhang, D. Probability Theory and Stochastic Processes; Science Publications: Beijing, China, 2005.
6. Doi, M. Second quantization representation for classical many-particle system. J. Phys. A: Math. Gen. 1976, 9, 1465. https: //doi.org/10.1088/0305-4470/9/9/008.
7. Trimper, S. Master equation and two heat reservoirs. Phys. Rev. E 2006, 74. https://doi.org/10.1103/physreve.74.051121. 352
8. Peliti, L. Path integral approach to birth-death processes on a lattice. J. Phys. France 1985, 46, 1469–1483. https://doi.org/10.105 353 1/jphys:019850046090146900.
9. Garcia-Palacios, J.L. Introduction to the Theory of Stochastic Processes and Brownian motion problems, 2007, [arXiv:cond-mat.stat-mech/cond-mat/0701242].
10. Bach, R.; Pope, D.; Liou, S.H.; Batelaan, H. Controlled double-slit electron diffraction. New J. Physics 2013, 15, 033018
11. Aharonov, Y.; Bohm, D. Significance of Electromagnetic Potentials in the Quantum Theory. Phys. Rev. 1959, 115, 485–491.
12. Kleinert, H. Path Integrals in Quantum Mechanics, Statistics, Polymer Physics, and Financial Markets, 5th ed.; WORLD SCIENTIFIC, 2009; [https://www.worldscientific.com/doi/pdf/10.1142/7305].
13. Berestycki, N.; Sousi, P. Applied Probability - Online Lecture Notes, 2007. http://www.statslab.cam.ac.uk/ ps422/notes-new.pdf.

14. Kosztin, I.; Faber, B.; Schulten, K. Introduction to the diffusion Monte Carlo method. Am. J. Phys. 1996, 64, 633–644. 360 https://doi.org/10.1119/1.18168. 361
15. Kaushal, R.S. Classical and quantum mechanics of complex Hamiltonian systems: An extended complex phase space approach. Pramana 2009, 73, 287–297. https://doi.org/10.1007/s12043-009-0120-x.
16. Rajeev, S.G. Dissipative Mechanics Using Complex-Valued Hamiltonians, 2007, [arXiv:quant-ph/quant-ph/0701141].
17. Morgan, J.D.; Simon, B. The calculation of molecular resonances by complex scaling. J. Phys. B At. Mol. Opt. Phys. 1981, 14, L167. https://doi.org/10.1088/0022-3700/14/5/002.
18. Botten, L.C.; Craig, M.S.; McPhedran, R.C.; Adams, J.L.; Andrewartha, J.R. The Finitely Conducting Lamellar Diffraction Grating. Optica Acta: International Journal of Optics 1981, 28, 1087–1102, [https://doi.org/10.1080/713820680].
19. Zloshchastiev, K.G. Quantum-statistical approach to electromagnetic wave propagation and dissipation inside dielectric media and nanophotonic and plasmonic waveguides. Phys. Rev. B 2016, 94. https://doi.org/10.1103/physrevb.94.115136.
20. Zloshchastiev, K.G. Generalization of the Schrödinger Equation for Open Systems Based on the Quantum-Statistical Approach. Universe 2024, 10. https://doi.org/10.3390/universe10010036.
21. Zloshchastiev, K. PROJECT Density Operator Approach for non-Hermitian Hamiltonians.
22. Sergi, A.; Zloshchastiev, K.G. Quantum entropy of systems described by non-Hermitian Hamiltonians. J. Stat. Mech. Theory Exp.2016, 033102.
23. Barreiro, J.T.; Müller, M.; Schindler, P.; Nigg, D.; Monz, T.; Chwalla, M.; Hennrich, M.; Roos, C.F.; Zoller, P.; Blatt, R. An open-system quantum simulator with trapped ions. Nature 2011, 470, 486–491
24. Zheng, C. Universal quantum simulation of single-qubit nonunitary operators using duality quantum algorithm. Sci. Rep. 2021, 11, 3960.
25. Fertik, M.B.; Scott, T.; Dignan, T. Identifying Information related to a Particular Entity from Electronic Sources, using Dimensional Reduction and Quantum Clustering, 2014. US Patent No. 8,744,197.
26. Wang, S.; Dignan, T.G. Thematic Clustering, 2014. US Patent No. 888, 665, 1 B1.
27. Scott, T.C.; Therani, M.; Wang, X.M. Data Clustering with Quantum Mechanics. Mathematics 2017, 5. https://doi.org/10.3390/math5010005.
28. Maignan, A.; Scott, T.C. A Comprehensive Analysis of Quantum Clustering: Finding All the Potential Minima. International Journal of Data Mining & Knowledge Management Process 2021, 11, 33–54.
29. Kumar, D.; Scott, T.C.; Quraishy, A.; Kashif, S.M.; Qadeer, R.; Anum, G. The Gender-Oriented Perspective in the Development of Type 2 Diabetes Mellitus Complications; SGOP Study. JPTCP 2023, 30, 2308–2318.